\documentclass{amsart}
\usepackage{graphicx}
\usepackage{amssymb}
\newtheorem{thm}{Theorem}[section]
\newtheorem{cor}{Corollary}[section]

\newtheorem{prop}{Proposition}[section]
\theoremstyle{definition}
\newtheorem{defn}{Definition}[section]
\theoremstyle{remark}
\newtheorem{rem}{Remark}[section]
\numberwithin{equation}{section}

\begin{document}

\title{The symmetry axioms in Laguerre planes}%
\author{Jaros\l aw Kosiorek and Andrzej Matra\'{s}}%
\address{Jaros\l aw Kosiorek and Andrzej Matra\'{s}\\
Department of Mathematics and Informatics\\
UWM Olsztyn\\
\.Zo\l nierska 14\\
10-561 Olsztyn\\
Poland}%
\email{matras@uwm.edu.pl}


\begin{abstract}
We introduce two axioms in Laguerre geometry and prove that they
provide a characterization of miquelian planes over fields of the
characteristic different from 2. They allow to describe an
involutory automorphism that sheds some new light on a Laguerre
inversion as well as on a symmetry with respect to a pair of
generators.
\end{abstract}
\maketitle

\section*{Introduction}
Classical Benz planes over commutative fields are described by Miquel's Theorem.
For Minkowski planes the simplification of this is given
by the symmetry postulate (S) which defines a natural orthogonal relation
and guarantees strong properties of the symmetry with respect to a circle \cite{HK1}.

In \cite{A1} R. Artzy has introduced the symmetry postulate ($\Pi$) for Laguerre planes.
The relation between (S) and ($\Pi$) can be seen through
the properties of ovals associated with circles of the plane .

But \cite{A1} has not contain a description of an automorphism
realizing the symmetry of the configuration of ($\Pi$).

The postulate $(\Pi)$ suggests that the tangent relation can be used to the classification
of Laguerre planes and to define involutory automorphisms.
Developing this idea we introduce two
axioms describing miquelian Laguerre planes (of characteristic different from 2) in terms of
tangent relation.
The postulate like ours were considered in various circle geometries (cf. \cite{F}, \cite{T}).
However it was unknown
that they imply the Miquel's Theorem in the Laguerre case.
This approach gives possibility to the define an involutory automorphism
that can be considered as the generalization of
both the Laguerre symmetry with two pointwiese fixed generators (cf. \cite{H})
and the Laguerre inversion (cf. \cite{M}). We called it the double tangency symmetry.
This automorphism realizes the symmetry of ($\Pi$)
and even allows to strengthen it.

In the paper we consider pencils $\langle K,L\rangle$ of circles
commonly tangent to given non-tangent circles $K$ and $L$. They may
be treated as analogous to notion of pencils of circles orthogonal
to two circles of a Minkowski plane. In the case of a pythagorean
field the geometry in the set of circles fixed by double tangency
symmetry without fixed points is M\"{o}bius geometry that contains
pencils $\langle K,L\rangle$ as circles.

\section{Notations and basic definitions}

A  \it Laguerre plane \/ \rm is a structure
$\mathbb{L}=(\mathcal{P},\mathcal{C},\parallel)$,
where $\mathcal{P}$ is a set of \it points \/\rm denoted by small Latin letters,
$\mathcal{C}\subset 2^{\mathcal{P}}$ is a set of
\it circles \/\rm denoted by capital Latin letters,
and $\parallel\subset\mathcal{P}\times\mathcal{P}$ is an equivalence relation known as
\it parallelity \rm.
The equivalence classes of the relation $\parallel$ will be called \it
generators \/ \rm and denoted by capital Latin letters.
We suppose that the following axioms are satisfied:
\begin{enumerate}
\item[(1)] Any three mutually non-parallel points $a,b,c$ are joined by a
unique circle, we denote this circle by $(a,b,c)^{\circ}$.
\item[(2)]
For every circle $K$ and any two non-parallel points $p\in K, q\notin K$
there is precisely one circle
$L$ which passes through $q$ such that $K\cap L=\{p\}$.

\item[(3)]
For any point $p$ and and any circle $K$ there exists exactly one point $q$
such that $p\parallel q$ and $q\in K$, we denote it by $pK$.
\item[(4)]
There is a circle containing at least three but not all points.
\end{enumerate}
We say circles $K$ and $L$ are tangent at $p$ iff $K\cap L=\{p\}$ or $K=L$.
If $p$ is a point of a circle $K$ then
we write  $\langle p,K\rangle$ for the pencil of circles tangent to
$K$ at the point $p$. If $q\nparallel p$ the circle of the pencil
$\langle p,K\rangle$ passing through $q$ will be denoted by $(p,K,q)^{\circ}$.
For any pair of non-parallel points $x,y$
the set of circles containing them will be called the \it pencil of
circles with the vertexes \rm $x,y$ and denoted by $\langle x,y\rangle$.

\begin{defn} \label{d1.1}
An ordered quadruple of points $a,b,c,d$ is said to be {\it
concyclic} if $a,b,c,d\in K$ for a certain circle $K$ or
$a\parallel b$, $c\parallel d$ and $a \nparallel c$.
Such quadruples we denote by
$(a,b,c,d)_{\triangle}$ \cite{C}. If the points $a,b,c,d$
belong to one circle we say that $a,b,c,d$ are {\it properly
concyclic}.
\end{defn}

The class of Laguerre planes embedded in projective spaces over commutative fields
was well investigated and characterized by Miquel's Theorem \cite{C}.
\begin{enumerate}
\item[(M)]
For  any eight different points $a,b,c,d,e,f,g,h$
the relations
$(a,c,b,d)_{\triangle}$, $(a,e,b,h)_{\triangle}$,
$(a,g,d,h)_{\triangle}$, $(b,f,c,e)_{\triangle}$, $(c,g,d,f)_{\triangle}$
imply that $(e,g,f,h)_{\triangle}$.
\end{enumerate}
The class of ovoidal Laguerre planes embedded in projective spaces over skew-fields
was described by the Bundle Theorem \cite{K}.
\begin{enumerate}
\item[(B)]
For any eight different points $a,b,c,d,e,f,g,h$ from
$(a,c,b,d)_{\triangle}$, $(c,e,d,f)_{\triangle}$,
$(e,g,f,h)_{\triangle}$, $(g,a,h,b)_{\triangle}$, $(a,e,b,f)_{\triangle}$
follows that $(c,g,d,h)_{\triangle}$.
\end{enumerate}
We note that in the Bundle Theorem at most two of the six quadruples can be not properly concyclic.
We will use the following weak form of this axiom (Fig. 1). This axiom is known as the
degenerated form of Miquel's Theorem and denoted by (M2)in \cite{Br1}, \cite{Kn}. It is weaker then Bundle Theorem
and characterizes so called \it elation Laguerre planes \rm. It follows from this remark that
every ovoidal Laguerre plane is an elation plane. Other proof of this statement can be found in \cite{S1}.

\begin{figure}[!h]
\includegraphics[width=0.3\textwidth]{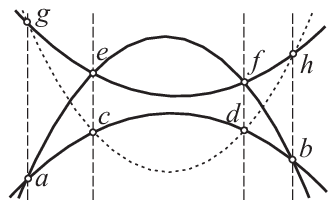}
\caption{} \label{r1}
\end{figure}

In \cite{A1} R. Artzy have proved that Miquel's Theorem is equivalent
to the following condition (Fig. 2):
\begin{enumerate}
\item[($\Pi$)]
Let $a,b,c,x$ be mutually non-parallel points such that
$x\notin(a,b,c)^{\circ}$ and let $p,q$ be points such that
$p\parallel c$, $p\in(a,b,x)_{\circ}$, $q\parallel b$, $q\in(a,c,x)_{\circ}$.
If $K$ is a circle containing $x$ tangent to $(a,b,c)^{\circ}$ in $a$,
then $K\cap(p,q,x)^{\circ}=\{x\}$.
\end{enumerate}

\begin{figure}[!h]
\includegraphics[width=0.5\textwidth]{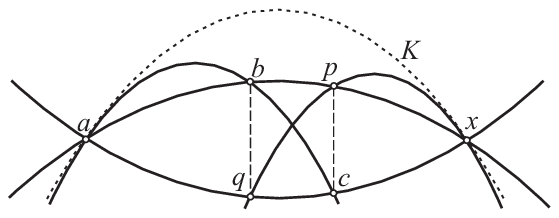}
\caption{} \label{r1}
\end{figure}


A miguelian Laguerre plane will be called {\it of characteristic 2}
if the corresponding field is of the same characteristic.

A bijective transformation of $\mathcal{P}$ preserving the class of circles
is said to be an {\it automorphism} of $\mathbb{L}$.
If $X,Y$ are different generators and $M$ is a circle
then an involutory automorphism fixing the generators $X,Y$ (pointwise)
and the circle $M$ (not pointwise) is called a {\it Laguerre symmetry} and
denoted by the symbol $\mathrm{S}_{X,Y;M}$.

H. M\"{a}urer \cite{M} constructed a Laguerre plane in
a M\"{o}bius plane $\mathbb{M}=(\mathcal{P},\mathcal{C},\in)$ over a pythagorean field.
Let us take a point $p\in \mathcal{P}$ and consider
the set $\mathcal{P}'$ of oriented circles passing through $p$,
and the set $\mathcal{C}'$ consisting of all oriented circles not containing $p$
and all points different from $p$.
A circle $A\in \mathcal{P}'$ is incident to $q\in \mathcal{C}'$ iff
$q\in A$, and $A$ is incident to a circle $B\in \mathcal{C}'$ iff
$A,B$ are tangent(with consistent orientation).
The same construction was applied for flat M\"{o}bius planes \cite{G1}.

In the Laguerre plane considered above
there exists a unique involutory automorphism without fixed points
such that all circles corresponding to the points of M\"{o}bius plane are preserved,
but all orientations are changed.
This automorphism is known as  {\it Laguerre inversion} \cite{M}.

\begin{figure}[!h]
\includegraphics[width=0.3\textwidth]{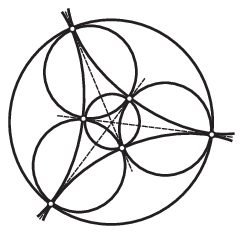}
\caption{} \label{r3}
\end{figure}

We will need the following.
\begin{prop} \label{p1.1}
Suppose that $\mathbb{L}$ is a miquelian Laguerre plane of characteristic $2$.
If $p\in K,L$ and $K,L$ are circles tangent to a circle $M$,
then $K,L$ are tangent in $p$.
\end{prop}

\begin{proof}
Let us consider the derived affine plane at $p$ and
its projective extension.
The circle $M$ induces a parabolic curve and circles $K,L$
induce lines tangent to this curve.
The lines of the extension induced by $K,L$ meet
the line in infinity at the same point
(this point is the kernel of the oval associated with $M$ cf. \cite{G3}).
This means that these lines are parallel in the affine plane, hence the
circles $K,L$ are tangent.
\end{proof}

Proposition \ref{p1.1} is illustrated on the following "tangential"
model of the Laguerre plane of order 2 (Fig. 3).

\section{Axioms and their representation theorem}
Let us consider the following axioms:
\begin{enumerate}
\item[(C)] For any circles $K,L$ and any point $p\in K\setminus L$ there exists exactly one circle $M$ such that $M\in
\langle p,K\rangle$ and $\overline{\overline{M\cap L}}=1$ (Fig. 4).
\begin{figure}[!h]
\includegraphics[width=0.4\textwidth]{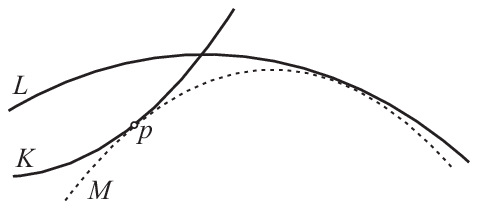}
\caption{} \label{r4}
\end{figure}
\item[(S)] If $K,L,M,N$ are circles and $a,b,c,d$ are points
such that $K\cap L=\{a\}$, $L\cap M=\{b\}$, $M\cap N=\{c\}$, $N\cap K=\{d\}$
and $a\nparallel c$, then there exists a circle passing through $a,b,c,d$ (Fig. 5).
\begin{figure}[!h]
\includegraphics[width=0.35\textwidth]{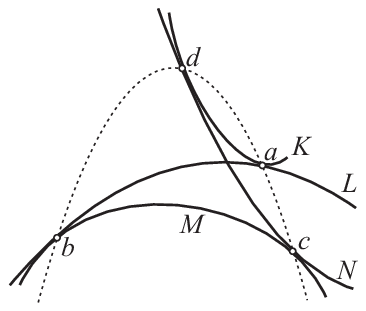}
\caption{} \label{r5}
\end{figure}
\end{enumerate}

The first axiom was named (C), because it does not hold for planes
of characteristic $2$. In the case when $K,L$ are not tangent, the
axiom (C) can be drawn from 4.4 of \cite{F} for miquelian Laguerre
planes of characteristic distinct from 2. The class of planes
satisfying (C) is more larger. Note that (C) is a special case of
the solution of Apollonius problem for connected topological
Laguerre planes of finite dimension (cf. \cite{S}, Theorem 7.1).
This statement distinguish Laguerre planes among analogical
classes of other Benz planes. For Minkowski planes and
M\"{o}bius planes (both of characteristic different form 2) if $K$ and $L$ are tangent (not in $p$)) there exist
exactly one $M$ (different from $K$) satisfying the assertion of
(C).

The second axiom we denote by (S) because as we show it provides
the symmetry and implies the symmetry axiom ($\Pi$) of Artzy (cf.
\cite{A1}). Similar condition was considered in \cite{T} for
M\"{o}bius planes, as one of the special forms of Miquel's Theorem.
We prove that, in Laguerre geometry, this degenerate version of
the Miquel's Theorem (together with (C)) implies the plane is miquelian.

The axiom (C) allows us to introduce for any distinct circles $K,L$ and a point
$p\in K\setminus L$ the following
notation:
$(p,K,L)^{\circ}$ - the unique circle of a pencil $\langle p,K\rangle$ tangent to $L$,
and $pKL$ - the point of tangency of circles $L$ and $(p,K,L)^{\circ}$.
For the case when $p\in K\cap L$ we define $pKL=p$.

If $K,L$ are any non-tangent circles we define a bijection $h^{L}_{K}:K\rightarrow L$ by the formula $h^{L}_{K}(x)=xKL$
for $x\in K$. By definition $h^{L}_{K}\circ h^{K}_{L}=id_{K}$. For distinct circles $K,L$ the set of circles tangent to
$K$ and $L$ we denote by $\langle K,L\rangle$ and call \it double tangency pencil. \rm\\
If we apply (C) for $K,L$ tangent (not in $p$) we obtain\\

\begin{prop} \label{p2.1}
Three mutually tangent circles are tangent at the same point.
\end{prop}

From Proposition \ref{p2.1} follows that if $K$ is tangent to $L$ at a point $p$,
then $\langle K,L\rangle=\langle p,K\rangle$.\\

\begin{prop} \label{p2.2}
If $K,L,M,N$ are circles
and $a,b,c,d$ are points
such that $K\cap L=\{a\}$, $L\cap M=\{b\}$, $M\cap N=\{c\}$, $N\cap K=\{d\}$ and $a\parallel c$, then $b\parallel d$
\end{prop}

\begin{proof}
If $b\nparallel d$, then the axiom (C) guarantees the existence of
a circle passing through $a,b,c,d$ which
contradicts to $a\parallel c$.
\end{proof}

As a consequence of (S) and Proposition \ref{p2.2} we obtain
\begin{cor} \label{c2.1}
If $K\cap L=\{a\}$, $L\cap M=\{b\}$, $M\cap N=\{c\}$, $N\cap K=\{d\}$, then $(a,c,b,d)_{\triangle}$.
\end{cor}
To prove the Artzy's postulate ($\Pi$) we need the following reformulation:
\begin{enumerate}
\item[($\Pi'$)]
Let $a,b,c,x$ are mutually non-parallel points such that $x\notin(a,b,c)^{\circ}$ and $q\parallel b$,
$q\in (a,c,x)^{\circ}$.
If $K$ is a circle passing through $x$ tangent to $(a,b,c)^{\circ}$ at $a$
and $L$ is a circle passing through $q$ tangent to $K$ at $x$,
then $L\cap(a,b,c)^{\circ}=\{p,x\}$ where $p\parallel c$.
\end{enumerate}

\begin{thm} \label{t2.1}
The axioms {\rm (C)} and {\rm (S)} implies {\rm($\Pi$)}.
\end{thm}

\begin{proof}
Let $a,b,c,x,q$ and $L,K$ be as in  the postulate ($\Pi'$).
According to (C)
there exists a circle $M$ tangent to $(a,b,c)^{\circ}$ at $b$ and tangent to $L$ in a certain point $p$.
It follows from (S) that $p\in (a,b,x)^{\circ}$.
The axiom (C) implies the existence of a circle $N$ tangent to $L$ at $q$
and tangent to the circle $(a,b,c)^{\circ}$ at a point $c'$.
By (S), $c'\in(a,q,x)^{\circ}$.
We have $c'\neq a$
(otherwise $(a,b,c)^{\circ}$ and $N$ are tangent at $a$,
and $N,K,L$ is a triple of circles mutually tangent in distinct points,
and the existence of such triple contradicts to Proposition \ref{p2.1}.
Hence $c=c'$ and Proposition \ref{p2.2} gives the claim.
\end{proof}

\begin{thm} \label{t2.2}
A Laguerre plane $\mathbb{L}$ satisfies the axioms $\mathrm{(C)}$ and $\mathrm{(S)}$ iff
it is a miquelian plane of characteristic different from 2.
\end{thm}

\begin{proof}
$\Rightarrow$
$\mathbb{L}$ is miquelian according to Theorem \ref{t2.1}
and  \cite{A1} (Theorem 5, p. 115).
By Proposition \ref{p2.1} and Proposition \ref{p1.1},
the characteristic of $\mathbb{L}$ is not $2$.

$\Leftarrow$
The statement can be obtained by an immediate verification
from the analytical description of miguelian Laguerre planes.
\end{proof}

Now we can draw something more for the configuration of ($\Pi$) (for an arbitrary characteristic).

\begin{thm} \label{t2.3}
Let $\mathcal{L}$ be miquelian and $a,b,c,x,p,q$ be as in {\rm ($\Pi$)}.
Then the circle tangent to $(q,p,x)^{\circ}$ at $p$ passing through $p$ is tangent to $(a,b,c)^{\circ}$ at $b$.
\end{thm}

\begin{proof}
For case when the characteristic is distinct from $2$, our statement follows from ($\Pi$) and Theorem \ref{t2.2}.
If the characteristic is $2$ then we use Proposition \ref{p1.1}.

Let $M$ be a circle passing through $x$ and tangent to $(a,b,c)^{\circ}$ at $b$.
Proposition \ref{p1.1} shows that
$M$ is tangent to $K$ at $x$,
hence it is tangent to $(q,p,x)^{\circ}$ at $x$.
Let $N$ be the circle tangent to
$(q,p,x)^{\circ}$ at $p$ and passing through $b$.
By Proposition \ref{p1.1}, $N$ is tangent to $M$ in $b$
(since $N,M$ are two circles tangent to $(q,p,x)^{\circ}$ passing through $b$).
Hence $N$ is tangent to $(a,b,c)^{\circ}$ at $b$.
\end{proof}

\section{Double tangency symmetry}
Throughout this section we suppose that $\mathbb{L}$ satisfies (C) and (S).
We define an involutory automorphism which generalizes both Laguerre symmetry (cf. \cite{H}) and
Laguerre inversion (cf. \cite{M}).

\begin{thm} \label{t3.1}
Let $K,L$ be non-tangent circles of $\mathbb{L}$.
There exists an involutory automorphism $\phi$ with properties
\begin{enumerate}
\item[(1)] $\phi(x)=xKL$ for all $x\in K$,
\item[(2)] if $\phi(x)\neq x$ and $M\in\langle x,\phi(x)\rangle$, then $\phi(M)=M$.
\end{enumerate}
\end{thm}

\begin{proof}
We define the automorphism $\phi$ by the formula (Fig. 6)
\begin{equation}
\phi(x)=((xK)KL)(x,y,yKL)^{\circ},
\end{equation}
where $y$ is an arbitrary point such that $y\in K\setminus L$,
$y\nparallel x$.

\begin{figure}[!h]
\includegraphics[width=0.45\textwidth]{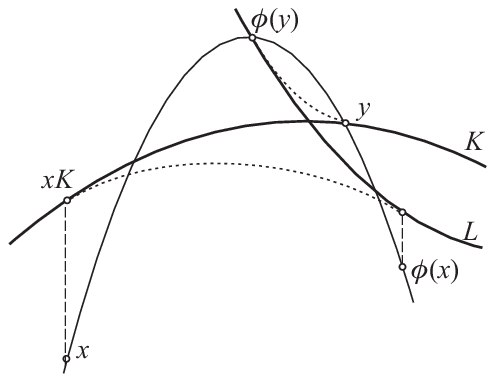}
\caption{} \label{r6}
\end{figure}

First we show that the definition does not depend from $y$.

By (S), there is a circle containing $x,y,xKL,yKL$;
this implies that $\phi(x)=xKL$ for every $x\in K$.
Similarly, we have $\phi(x)=xLK$ for every $x\in L$
($xLK\parallel(xK)KL$ according to Proposition \ref{p2.2}).

Now consider the case when $x\notin K\cup L$.
Let us take an arbitrary point $z\in K\setminus L$ such that $z\neq y$, $z\nparallel x$
and define
$$x':=\phi(x),\; x'':=(xK)KL,\; y':=yKL,\; z':=zKL.$$
By our construction, the quadruple
$$x,\;x',\;y,\;y'$$
is contained in a circle.
According to (S)
the same holds for the quadruples
$$x,\;x',\;xK,\;x'',$$
$$z,\;z',\;xK,\;x'',$$
$$z,\;z',\;y,\;y',$$
$$y,\;y',\;xK,\;x''.$$
By the bundle theorem (B), the points $z,z',x,x'$ are concyclic.
This shows that $\phi$ is well-defined for $x\notin K\cup L$.

The transformation $\phi$ is an involution; indeed, $(x'K)KL\parallel x$ and hence $\phi(x')=x$.

It follows from the definition of $\phi$ that $\phi(M)=M$ for each $M\in\langle y, yKL\rangle$, where $y\in K$.
Now we want to establish the same for all $y\notin K\cup L$.

Let $x\in M\in\langle y,\phi(y)\rangle$ and $x \neq y,\phi(y)$.
Now to define $\phi(x)$ and $\phi(y)$ we will use the points $yK$ and $xK$ (respectively).
We have (Fig. 7)
\begin{figure}[!h]
\includegraphics[width=0.8\textwidth]{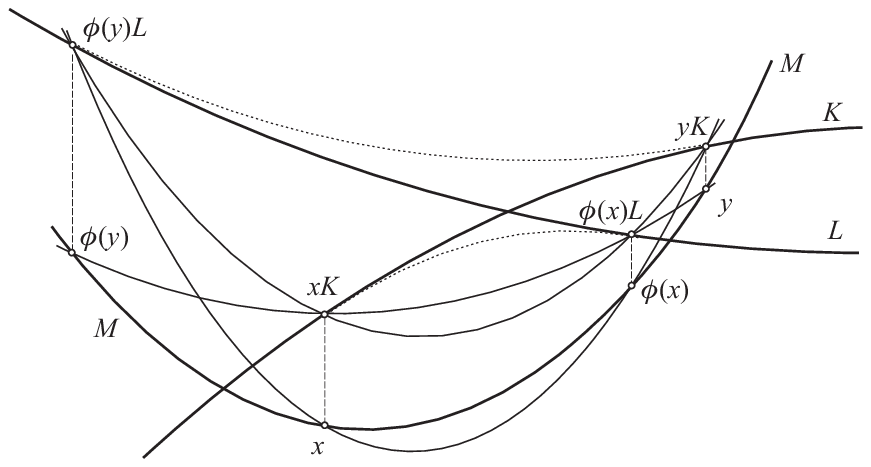}
\caption{} \label{r7}
\end{figure}
$$(y,yK,\phi(y),\phi(y)L)_{\triangle}\;\mbox{ (a pair of generators)},$$
$$(y,xK,\phi(x)L,\phi(y))_{\triangle}\;\mbox{ (a circle)},$$
$$(x,xK,\phi(x),\phi(x)L)_{\triangle}\;\mbox{ (a pair of generators)},$$
$$(yK,x,\phi(x),\phi(y)L)_{\triangle}\;\mbox{ (a circle)}.$$
By (S), we have $(yK,xK,\phi(y)L,\phi(y)L)_{\triangle}$ and from Bundle Theorem
we draw that $(y,x,\phi(y),\phi(x)$ are properly concyclic and the circle $M$ is invariant for $\phi$.

We show now that $\phi$ is an automorphism. Let $x,y,z,t$ be properly concyclic points.
By the arguments  given above,
the points of the following quadruples are properly concyclic:
$$x,\;y,\;\phi(x),\;\phi(y),$$
$$y,\;z,\;\phi(y),\;\phi(z),$$
$$z,\;t,\;\phi(z),\;\phi(t),$$
$$t,\;x,\;\phi(t),\;\phi(x).$$
According to Miquel's Theorem, $\phi(x),\phi(y),\phi(z),\phi(t)$ are properly concyclic.

In the end, we prove that there is only one automorphism
satisfying the conditions of the theorem.
This follows directly from the following:
\begin{equation}\label{eq-2}
\phi(x)=xM\phi(M)
\end{equation}
for any circle $M$ such that $M\neq\phi(M)$ and any $x\in M$.

Suppose first $x\notin\phi(M)$, so $x\neq\phi(x)$. Let us prove \eqref{eq-2}. We consider
the circle $N:=(x,M,\phi(x))^{\circ}$. The condition (2)
guarantees that $\phi$ preserves $N$ because $N\in\langle x,\phi(x)\rangle$.
Hence $N$ is tangent to $\phi(M)$ at the point
$\phi(x)$. Similarly we draw that the tangency o $M$ and $\phi(M)$ is impossible
(contradiction with Proposition\ref{p2.1}).
If $x\in M\cap\phi(M)$, then $\overline{\overline{M\cap\phi(M)}}=2$ and
 we obtain $\phi(x)=x$ because in the other case $M,\phi(M)\in\langle x,\phi(x)\rangle$, hence $M=\phi(M)$.
\end{proof}

The automorphism defined in Theorem \ref{t3.1}
will be called {\it the double tangency symmetry associated with $K,L$}
and denoted by $\mathrm{S}_{\{K,L\}}$.

\begin{thm} \label{t3.2}
An automorphism $\sigma$ of $\mathbb{L}$ is a Laguerre symmetry  if and only if
$\sigma=\mathrm{S}_{\{K,L\}}$, where
$|K\cap L|=2$.
\end{thm}

\begin{proof}
($\Leftarrow$).
Let $K,L$ be circles such that $K\cap L=\{p,q\}$.
We take an arbitrary $x\parallel p$.
From the definition of $\mathrm{S}_{\{K,L\}}$ we obtain $(xK)KL=p$,
hence
$$\mathrm{S}_{\{K,L\}}(x)=p(x,y,yKL)^{\circ}=x.$$
Similarly we conclude $\mathrm{S}_{\{K,L\}}=x$ for $x\parallel q$.
Then $\mathrm{S}_{\{K,L\}}$ coincides with $\mathrm{S}_{P,Q;M}$
where $P,Q$ are the generators containing $p,q$ (respectively),
and $M$ is an arbitrary circle preserved by $\mathrm{S}_{\{K,L\}}$.

($\Rightarrow$).
Let $\mathrm{S}_{P,Q;M}$ be a Laguerre symmetry.
There is  a circle $K$ such that
$$K\cap(P\cup Q)=M\cap(P\cup Q),$$
where $K\neq M$.
We define $L:=\mathrm{S}_{P,Q;M}(K)$;
it is trivial that $L\neq K$. Let us take $x\in K\setminus\{P\cup Q\}$. We have
$$\mathrm{S}_{P,Q;M}(x)=x'\in L.$$
The circle $N:=(x,K,x')^{\circ}$ is invariant by $\mathrm{S}_{P,Q;M}$, because $N\in\langle x,x'\rangle$.
Thus  $N$ is tangent to $L$ at $x'$, so
$x'=\mathrm{S}_{\{K,L\}}(x)$.
By the part ($\Leftarrow$),
$\mathrm{S}_{\{K,L\}}(x)=\mathrm{S}_{P,Q;M}(x)$ for $x\in P\cup Q$.
Since
$\mathrm{S}_{\{K,L\}}(x)=\mathrm{S}_{P,Q;M}(x)$ for $x\in P\cup Q\cup K\cup L$, we get
$\mathrm{S}_{\{K,L\}}=\mathrm{S}_{P,Q;M}$
\end{proof}

\begin{rem} \label{r3.1}
The symmetries associated with  the configuration of Artzy's axiom ($\Pi$) are always
$\mathrm{S}_{\{K,L\}}$, where
$K=(a,b,c)^{\circ}$, $L=(x,p,q)^{\circ}$.
A symmetry with respect to a pair of generators explanes the configuration
only for the case when $K$ and $L$ are intersecting.
\end{rem}

\begin{rem}\label{r3.2}
In Section 1 we have considered an example of a Laguerre plane
constructed in M\"{o}bius plane.
For this case, the Laguerre inversion $\mathrm{S}_{\{K,L\}}$ with $K\cap L=\emptyset$ have
the following attractive interpretation.
Suppose that $K,L$ are the
circles corresponding to a circle $M$ ($p\notin M$) of M\"{o}bius plane with opposite orientations.
We get the following:
\begin{enumerate}
\item[---]
$\langle K,L\rangle$ corresponds to the set of points of the circle $M$,
\item[---]
all fixed circles of $\mathrm{S}_{\{K,L\}}$ corrrespond to points of the M\"{o}bius plane,
\item[---]
any $\langle K',L'\rangle$ with $\mathrm{S}_{\{K,L\}}(K')=L'$ corresponds to the set of points of a circle not
containing $p$,
\item[---]
any pencil $\langle x,y\rangle$ (where $\mathrm{S}_{\{K,L\}}(x)=y$) corresponds to the set of points
of a circle passing through $p$,
\item[---]
$\mathrm{S}_{\{K,L\}}$ is the Laguerre inversion.
\end{enumerate}
Hence the geometry of circles fixed by $\mathrm{S}_{\{K,L\}}$ in this case can be seen as the M\"{o}bius geometry if we
take as M\"{o}bius circles double tangency pencils $\langle M,N\rangle$ and pencils $\langle x,y\rangle$
completed by a point $\{\infty\}$.
\end{rem}

\noindent
{\bf Acknowledgement}.
The authors thank Mark Pankov for  useful discussions and valuable comments.


\end{document}